\documentclass{article}%
\usepackage{amsmath}%
\usepackage{amsfonts}%
\usepackage{amssymb}%
\usepackage{graphicx}
\begin{document}

\title{Derivatives with respect to the order of the Legendre Polynomials}
\author{Bernard J. Laurenzi\\Department of Chemistry, UAlbany, The State University of New York \\1400 Washington Ave., Albany, N. Y. 12222}
\date{December 28, 2014}
\maketitle

\begin{abstract}
Expressions for the derivatives of the Legendre polynomials of the first kind
with respect to the order of these polynomials i.e. $\mathsf{P}_{n}%
(z)=[\partial^{\,n}P_{\nu}(z)/\partial\,\nu^{\,n}]_{\nu=0}$ \ are given. \ An
explicit form for the fourth derivative is presented.

\end{abstract}

\section{Introduction}

Recently, R. Szmytkowski \cite{RSZM} has obtained expressions for the first
three derivatives of the Legendre functions of the first kind $P_{\nu}(z)$
with respect to order $\nu$ i.e. $\mathsf{P}_{n}(z)=[\partial^{\,n}P_{\nu
}(z)/\partial\,\nu^{\,n}]_{\nu=0}\;$for $1\leq\nu\leq3.$ \ That is to say,%
\begin{align*}
\mathsf{P}_{0}(z)  & =1,\\
\mathsf{P}_{1}(z)  & =\ln(\tfrac{1+z}{2})=-\,Li_{1}(\tfrac{1\,-\,z)}{2}),\\
\mathsf{P}_{2}(z)  & =-\,2\,Li_{2}(\tfrac{1-z}{2}),\\
\mathsf{P}_{3}(z)  & =12\,Li_{3}(\tfrac{1+z}{2})-6\ln(\tfrac{1+z}{2}%
)\,Li_{2}(\tfrac{1+z}{2})-\pi^{2}\ln(\tfrac{1+z}{2})-12\,\zeta(3),
\end{align*}
where $Li_{\mu}(z)$ is the polylogarithm function \cite{LEW} of order $\mu$
and $\zeta(s)$ is the Riemann zeta function. \ These derivatives arise in
studies of tidal hydrodynamics and are part of a recent and continuing
interest in the variation of well known polynomials and other higher
transcendental functions with respect to their orders \cite{BRYC}%
,\cite{SESM},\cite{COHL},\cite{BUSC},\cite{LAUR}.

\bigskip

In this note we will give a general expression for the derivatives
$\mathsf{P}_{n}(z)$ and discuss the expected increasing complexity of these
expressions as $n$ increases. \ 

\bigskip

The Legendre functions of the first kind $P_{\nu}(z)$ satisfy the differential
equation \cite{LEG}%
\[
\lbrack\frac{d}{d\,z}(1-z^{2})\frac{d}{d\,z}+\nu(\nu+1)]\,P_{\nu}(z)=0\text{ ,
\hspace{0.35in}}-1\leq z\leq1.
\]
Differentiating the latter expression with respect to $\nu$ and evaluating the
result at $\nu=0$ we get for $\mathsf{P}_{n}(z)$ the relation%
\[
\frac{d}{dz}\left[  (1-z^{2})\frac{d\,\,\mathsf{P}_{n}(z)}{d\,z}\right]
=-\,n\,\mathsf{P}_{n-1}(z)-n\,(n-1)\,\mathsf{P}_{n-2}(z).
\]
Expressions for the desired derivatives $\mathsf{P}_{n}(z)$ can then be
reduced to quadratures by%
\begin{equation}
\mathsf{P}_{n}(z)=-n\,\int\frac{d\,z}{1-z^{2}}\int^{z}[\mathsf{P}%
_{n-1}(z^{\prime})+\,(n-1)\,\mathsf{P}_{n-2}(z^{\prime})]\,dz^{\prime}%
+C\ln(\tfrac{1+z}{1-z}),\label{eq1}%
\end{equation}
where $C$ is a constant of integration. \ The inner or \textit{first}
integrals in (1) will be discussed further below. \ With the use of the
identity
\[
P_{\nu}(1)=1,
\]
from which it follows that for $n\geq1$
\[
\mathsf{P}_{n}(1)=0,
\]
the constants of integration which arise in (1) can be evaluated.

\section{The Expression for $\mathsf{P}_{4}(z)$}

\bigskip

In this case we have to evaluate the expression%
\[
\mathsf{P}_{4}(z)=C\ln(\tfrac{1+z}{1-z})-4\int\frac{d\,z}{1-z^{2}}\int
^{z}[\mathsf{P}_{3}(z^{\prime})+\,3\,\mathsf{P}_{2}(z^{\prime})]\,dz^{\prime
},
\]
or more explicitly
\[
\mathsf{P}_{4}(z)=C\ln(\tfrac{1+z}{1-z})-4\int\frac{I(z)}{1-z^{2}}\,d\,z,
\]
where
\begin{align*}
I(z)  & =\int^{z}[12\,Li_{3}(\tfrac{1+z^{\prime}}{2})-6\ln(\tfrac{1+z^{\prime
}}{2})\,Li_{2}(\tfrac{1+z^{\prime}}{2})-\pi^{2}\ln(\tfrac{1+z^{\prime}}%
{2})-12\,\zeta(3)-6\,Li_{2}(\tfrac{1-z^{\prime}}{2})]\,dz^{\prime}.\\
& .
\end{align*}
The integrals in $I(z)$ are well known and we find the somewhat remarkable
result%
\[
I(z)=(z+1)[12\,Li_{3}(\tfrac{1+z}{2})-6\ln(\tfrac{1+z}{2})\,Li_{2}(\tfrac
{1+z}{2})-\pi^{2}\ln(\tfrac{1+z}{2})-12\,\zeta(3)].
\]
In the outer i.e. \textit{second} integration we have for $\mathsf{P}_{4}(z)$
\[
\mathsf{P}_{4}(z)=C\ln(\tfrac{1+z}{1-z})-4\int\frac{d\,z}{1-z}[12\,Li_{3}%
(\tfrac{1+z}{2})-6\ln(\tfrac{1+z}{2})\,Li_{2}(\tfrac{1+z}{2})-\pi^{2}%
\ln(\tfrac{1+z}{2})-12\,\zeta(3)],
\]
where $C$ is a constant of integration. \ All but one of the integrals are
elementary. \ With a change of variable
\[
t=\frac{1+z}{2},\hspace{0.2in}0\leq t\leq1
\]
we get for $\mathsf{P}_{4}(z)$
\begin{equation}
\mathsf{P}_{4}(z)=C\ln(\tfrac{t}{1-t})+C^{\prime}+24\left[  Li_{2}(t)^{2}%
+2\ln(1-t)\left[  Li_{3}(t)-\zeta(3)\,\right]  +\tfrac{\pi^{2}}{6}%
Li_{2}(1-t)+\,\mathfrak{I}(t)\right]  ,\label{eq2}%
\end{equation}
where
\[
\mathfrak{I}(t)=\int\frac{\ln(t)\,Li_{2}(t)}{1-t}\,dt,
\]
and $C^{\prime}$ is the second constant of integration.\medskip

The integral $\mathfrak{I}(t)$ has been obtained within \textit{Mathematica}
\cite{MATH}. \ After some simplication we have%
\begin{align}
\mathfrak{I}(t)  & =[\ln^{2}(t)-\ln(t)\ln(1-t)]\,Li_{2}(t)-\ln^{2}%
(1-t)\,Li_{2}(1-t)+\ln^{2}(\tfrac{t}{1-t})\,Li_{2}(\tfrac{t}{t-1})-\tfrac
{1}{2}\,Li_{2}(t)^{2}\label{eq3}\\
& +2\ln(t)\,Li_{3}(t)+2\ln(1-t)\,Li_{3}(1-t)-2\ln(\tfrac{t}{1-t}%
)\,Li_{3}(\tfrac{t}{t-1})\nonumber\\
& +2\,\left[  Li_{4}(t)-\,Li_{4}(1-t)+\,Li_{4}(\tfrac{t}{t-1})\right]
+\ln^{2}(1-t)[\tfrac{1}{2}\ln^{2}(t)-\ln(t)\ln(1-t)+\tfrac{1}{4}\ln
^{2}(1-t)].\nonumber
\end{align}
Using the following identities for the di and trilogarithm functions%
\begin{align*}
Li_{2}(1-z)  & =-\,Li_{2}(z)+\pi^{2}/6-\ln(z)\ln(1-z),\\
Li_{2}(\tfrac{z}{z-1})  & =-\,Li_{2}(z)-\tfrac{1}{2}\ln^{2}(1-z),\\
Li_{3}(\tfrac{z}{z-1})+Li_{3}(1-z)+Li_{3}(z)-\zeta(3)  & =\tfrac{\pi^{2}}%
{6}\ln(1-z)-\tfrac{1}{2}\ln(z)\ln^{2}(1-z)+\tfrac{1}{6}\ln^{3}(1-z),
\end{align*}
and gathering terms in equations 2 and 3 we get\medskip\
\begin{equation}
\mathsf{P}_{4}(z)=\pi^{4}/15+24\left[
\begin{array}
[c]{c}%
\tfrac{1}{2}Li_{2}(t)^{2}-\tfrac{\pi^{2}}{6}Li_{2}(t)+2\,Li_{4}(t)-2\,Li_{4}%
(1-t)+2\ln(t)\,{\LARGE \{}Li_{3}(1-t)-\zeta(3){\LARGE \}}\\
+\frac{1}{12}\ln^{4}(1-t)+\tfrac{\pi^{2}}{6}\ln^{2}(1-t)+2\,Li_{4}%
(\frac{t}{t-1})\\
+\ln(t)\ln(1-t){\LARGE \{}Li_{2}(t)-\frac{\pi^{2}}{2}-\frac{1}{3}\ln
^{2}(1-t)+\ln(t)\ln(1-t){\LARGE \}}%
\end{array}
\right]  .\nonumber
\end{equation}
\bigskip

The constants $C$ and $C^{\prime}$ in equation 4 having been \ evaluated by
setting \ $\mathsf{P}_{4}(1)=0$ \ with the result that $C=0$ and $C^{\prime
}=\pi^{4}/15.$ We note that the expression in equation 4 contains the
complicated term $Li_{4}(\tfrac{t}{t-1}).$ \ In contrast to the corresponding
expressions for the polylogarithms $Li_{2}$ and $Li_{3}$ with the same
argument, this term does not appear to be able to be rewritten in terms of the
polylogarithm $Li_{4}$ with simpler arguments as noted by Lewin \cite{LEW2}.
\ As a consequence it does not bode well for any likelihood of obtaining
explicit analytic expressions for $\mathsf{P}_{n}(z)$ for $n\geq5$.\ 

\begin{center}
\appendix \textbf{Appendix A}
\end{center}

\bigskip

Each of the \textit{first} integrals i.e.%
\[
\mathcal{I}_{\eta}\mathcal{(}z\mathcal{)=}\int^{z}\mathsf{P}_{\eta}(z^{\prime
})\,d\,z^{\prime},
\]
appearing in equation 1 occur twice i.e. in successive calculations of the
quantities $\mathsf{P}_{n}(z)$ and $\mathsf{P}_{n+1}(z)$. \ The first few of
these quantities are given below. \ It should be noted however, that these
expressions are more complicated than the combinations
\[
\int^{z}[\,\mathsf{P}_{n-1}(z^{\prime})\,+\,(n-1)\,\mathsf{P}_{n-2}(z^{\prime
})\,]\,d\,z^{\prime}.
\]
This is due to an internal cancellation of terms within the combinations.
\ The latter circumstance may indicate a deeper issue involving the
$\ \mathcal{I}_{\eta}\mathcal{(}z\mathcal{)}$ \ integrals. We have
\[
\mathcal{I}_{1}\mathcal{(}z\mathcal{)=(}1+z\mathcal{)[}\ln\mathcal{(}%
\tfrac{1+z}{2})-1\mathcal{]},
\]%
\[
\mathcal{I}_{2}\mathcal{(}z\mathcal{)}=-2\,(1+z)[\ln(\tfrac{1+z}%
{2})-1]+2\,(1-z)\,Li_{2}(\tfrac{1-z}{2}),
\]%
\begin{align*}
\mathcal{I}_{3}\mathcal{(}z\mathcal{)}  & =6\,(1+z){\LARGE [}2\,Li_{3}%
(\tfrac{1+z}{2})+\tfrac{\pi^{2}}{6}-1+2\,\zeta(3)-\{Li_{2}(\tfrac{1+z}%
{2})+\tfrac{\pi^{2}}{6}-1\}\ln(\tfrac{1+z}{2}){\LARGE ]}\\
& +6\,(1-z)\,{\LARGE [}Li(\tfrac{1+z}{2})+\ln(\tfrac{1-z}{2})\ln(\tfrac
{1+z}{2}){\Large ],}%
\end{align*}
The expression for $\mathcal{I}_{4}\mathcal{(}z\mathcal{)}$ has been computed
using \textit{Mathematica} and contains in a condensed form seventy-four terms
including a fifth order polylogarithm function. \ That quantity will not be
displayed here for the sake of brevity.\ 

Below we include integrals which occur in the calculation of $\mathcal{I}%
_{4}\mathcal{(}z\mathcal{)}$ but are not directly available within
\textit{Mathematica} i.e.%

\begin{align*}
\int Li_{4}(\tfrac{t}{t-1})\,dt  & =-\tfrac{1}{2}Li_{2}(\tfrac{t}{t-1}%
)^{2}+t\,Li_{4}(\tfrac{t}{t-1})+\ln(1-t)\,Li_{3}(\tfrac{t}{t-1}),\\
\int Li_{2}(t)^{2}\,dt  & =-2+6t+6[1-t-\tfrac{\pi^{2}}{9}]\ln(1-t)-2[1-t-\ln
(t)]\ln^{2}(1-t)\\
& -2[t-(1+t)\ln(1-t)]\,Li_{2}(t)+t\,Li_{2}(t)^{2}+4Li_{3}(1-t),
\end{align*}%
\begin{align*}
& \int\ln^{2}(t)\ln^{2}(1-t)\,dt=\\
& -4+24x+12\,[1-x]\ln(1-x)-2\,[1-x]\ln(1-x)^{2}-\tfrac{1}{2}\ln(1-x)^{4}\\
& -12x\ln(x)-4\,[1-2x]\ln(1-x)\ln(x)-2x\ln(1-x)^{2}\ln(x)+2\ln(1-x)^{3}%
\ln(x)\\
& +[2-\ln(1-x)^{2}]\ln(x)^{2}-(1-x)[2-2\ln(1-x)+\ln(1-x)^{2}]\ln(x)^{2}\\
& +[4-4\ln(1-x)+2\ln(1-x)^{2}]\,Li_{2}(1-x)-[4-4\ln(x)+2\ln(x)^{2}%
]\,Li_{2}(x)\\
& -[2\ln(1-x)^{2}-4\ln(1-x)\ln(x)+2\ln(x)^{2}]\,Li_{2}(\tfrac{x}{x-1})\\
& -4\,[1-\ln(x)]\,Li_{3}(x)+4\,[1-\ln(1-x)]\,Li_{3}(1-x)+4\,[\ln
(x)-\ln(1-x)]\,Li_{3}(\tfrac{x}{x-1})\\
& +4\,Li_{4}(1-x)-4\,Li_{4}(x)-4\,Li_{4}(\tfrac{x}{x-1}).
\end{align*}

\begin{center}
\ \ \ \appendix \textbf{Appendix B}\ \medskip
\end{center}

The integral $\mathfrak{I}(z)$ is interesting in that it provides a way to
obtain a closed form expression for the slowly converging infinite sum
$\sum_{k=1}^{\infty}\frac{\Psi^{^{\prime}}(k)}{k^{2}}$ where $\Psi^{^{\prime}%
}(k)$ is the trigamma function \ \cite{TRI}. This sum does not appear to have
been previously reported in the literature and is included here. \ Using the
limiting values for $\mathfrak{I}(z)$ i.e.
\begin{align*}
\mathfrak{I}(1)  & =-\,\frac{11}{360}\pi^{4},\\
\mathfrak{I}(0)  & =-\,\frac{1}{45}\pi^{4},
\end{align*}
and the infinite series representation for the dilogarithm function which
occurs in $\mathfrak{I}(z)$ we get
\begin{align*}
\mathfrak{I}(1)-\mathfrak{I}(0)  & =\sum_{k=1}^{\infty}\frac{1}{k^{2}}\int
_{0}^{1}\frac{z^{k}\,\ln(z)}{1-z}dz,\\
-\frac{\pi^{4}}{120}  & =-\sum_{k=1}^{\infty}\frac{\Psi^{^{\prime}}%
(k+1)}{k^{2}}.
\end{align*}
Expanding the trigamma function with the use of the recurrence relation for
$\Psi^{^{\prime}}$ i.e.
\[
\Psi^{^{\prime}}(k+1)=\Psi^{^{\prime}}(k)-1/k^{2},
\]
together with the well known sum $%
%TCIMACRO{\tsum _{k=1}^{\infty}}%
%BeginExpansion
{\textstyle\sum_{k=1}^{\infty}}
%EndExpansion
1/k^{4}=\pi^{4}/90$ \ we get the value of the desired summation i.e.
\[
\sum_{k=1}^{\infty}\frac{\Psi^{^{\prime}}(k)}{k^{2}}=\frac{7}{360}\pi^{4}.
\]

\end{document}